\documentclass[a4paper,10pt]{article}
\usepackage{epsfig,graphicx,verbatim,amsmath,amsfonts,amssymb,amsthm,alltt,algorithm,algorithmic,multirow}
\usepackage[margin=1in,nohead]{geometry}
\usepackage{adjustbox}
\DeclareGraphicsExtensions{.png,.pdf,.jpg}
\begin{document}
\newcommand{\raf}[1]{(\ref{#1})}
\newcommand{\sch}{\choose }
\newcommand{\ul}[1]{\underline{#1}}
\newcommand{\comp}{\hbox{$<\kern -3pt >$}}
\newcommand{\ncomp}
		{\;\hbox{\hbox{/}\kern -9.5pt \hbox{$<\kern -3pt >$}}}
\newcommand{\mpt}{\partial}

\newcommand{\meet}
	       {\hbox{$\wedge \kern -5.75pt \raise 1.5pt \hbox{$.$}\,$}}
\newcommand{\Meet}
	     {\hbox{$\bigwedge \kern -8pt \raise 0.75pt \hbox{$.$}\:$}}
\newcommand{\pplus}
	    {\hbox{$+ \kern -6pt \raise 9.05pt \hbox{$.$}\:$}}
\newcommand{\ld}
	       {\hbox{$< \kern -6pt \raise 2pt \hbox{$.$}\,$}}

\newcommand{\bve}{\bigwedge}
\newcommand{\fml}{((\mu-L)/\gs)}
\newcommand{\fmgl}{((\mu-L)/(\ga\mu))}

\newcommand{\ctnt}{\hbox{$\,\hat{\ }\,$}}

\newcommand{\ct}{\centerline}
\newcommand{\sat}{\models}
\newcommand{\prvs}{\vdash}
\newcommand{\frc}{\hbox{$\parallel \kern -5.7pt \hbox{$-$}$}}
\renewcommand{\iff}{\leftrightarrow}
\newcommand{\ra}{\rightarrow}

\newcommand{\lharp}{\leftharpoonup}
\newcommand{\rharp}{\rightharpoonup}

\newcommand{\nfrc}{\not \kern -5pt \frc}
\newcommand{\rest}{\vbox{\hbox{$\:\kern -2pt\mathbin{\vert\kern-3.1pt\lower-1pt
   \hbox{$\mathsurround=0pt\mathchar"0012$}\kern-4pt}\:$}}}
\newcommand{\drest}{\rest\rest}

\newcommand{\rs}{\!\!\!\!\!}
\newcommand{\pr}{\prec \!\!\!)}
\newcommand{\suc}{\succ \!\!(}
\newcommand{\nin}{\in\!\!\!\!\!/}
\newcommand{\rang}{<\!\!\! )}
\newcommand{\lang}{>\!\!\!\!\! (}

\newcommand{\empt}{\emptyset}

\newcommand{\cntd}{\subseteq}
\newcommand{\cnts}{\supseteq}
\newcommand{\pcntda}{\lower5pt\hbox{$\stackrel{\subset}{\neq}$}}
\newcommand{\pcntdb}{\lower5pt\hbox{$\stackrel{\supset}{\neq}$}}

\newcommand{\pexp}[1]{\hbox{$ #1 $}}

\newcommand{\real}[1]{\hbox{\rm #1}}
\newcommand{\Cal}[1]{{\cal #1}}

\newcommand{\ml}{\ell}

\newcommand{\ga}{\alpha}
\newcommand{\gga}{\gamma}
\newcommand{\gb}{\beta}
\newcommand{\gd}{\delta}
\newcommand{\gk}{\kappa}
\newcommand{\get}{\eta}
\newcommand{\gep}{\varepsilon}
\newcommand{\gvp}{\varphi}
\newcommand{\gl}{\lambda}
\newcommand{\gL}{\Lambda}
\newcommand{\Gl}{\Lambda}
\newcommand{\gch}{\chi}
\newcommand{\gp}{\pi}
\newcommand{\gps}{\psi}
\newcommand{\gs}{\sigma}
\newcommand{\gr}{\rho}
\newcommand{\ges}{\varsigma}
\newcommand{\gS}{\Sigma}
\newcommand{\gt}{\theta}
\newcommand{\gT}{\Theta}
\newcommand{\go}{\omega}
\newcommand{\gO}{\Omega}
\newcommand{\gG}{\Gamma}
\newcommand{\gx}{\xi}
\newcommand{\gz}{\zeta}
\newcommand{\e}{{\rm e}}
\newcommand{\m}{{\rm m}}

\newcommand{\bA}{\hbox{\bf A}}
\newcommand{\bC}{\hbox{\bf C}}
\newcommand{\bL}{\hbox{\bf L}}
\newcommand{\bI}{\hbox{\bf I}}
\newcommand{\bR}{\hbox{\bf R}}
\newcommand{\bE}{\hbox{\bf E}}
\newcommand{\bB}{\hbox{\bf B}}
\newcommand{\bN}{\hbox{\bf N}}
\newcommand{\bZ}{\hbox{\bf Z}}
\newcommand{\bT}{\hbox{\bf T}}
\newcommand{\bQ}{\hbox{\bf Q}}
\newcommand{\bU}{\hbox{\bf U}}
\newcommand{\bu}{\hbox{\bf u}}
\newcommand{\bl}{\hbox{\bf l}}
\newcommand{\bt}{\hbox{\bf t}}
\newcommand{\bh}{\hbox{\bf h}}
\newcommand{\bH}{\hbox{\bf H}}
\newcommand{\bP}{\hbox{\bf P}}
\newcommand{\bS}{\hbox{\bf S}}
\newcommand{\bq}{\hbox{\bf q}}

\newcommand{\HC}{\hbox{\rm\bf H$^{\infty}+$c}} 

\newcommand{\CA}{{\cal A}}
\newcommand{\CB}{{\cal B}}
\newcommand{\CC}{{\cal C}}
\newcommand{\CD}{{\cal D}}
\newcommand{\CE}{{\cal E}}
\newcommand{\CF}{{\cal F}}
\newcommand{\CG}{{\cal G}}
\newcommand{\CH}{{\cal H}}
\newcommand{\CI}{{\cal I}}
\newcommand{\CJ}{{\cal J}}
\newcommand{\CK}{{\cal K}}
\newcommand{\CL}{{\cal L}}
\newcommand{\CM}{{\cal M}}
\newcommand{\CN}{{\cal N}}
\newcommand{\CO}{{\cal O}}
\newcommand{\CP}{{\cal P}}
\newcommand{\CQ}{{\cal Q}}
\newcommand{\CR}{{\cal R}}
\newcommand{\CS}{{\cal S}}
\newcommand{\CT}{{\cal T}}
\newcommand{\CU}{{\cal U}}
\newcommand{\CV}{{\cal V}}
\newcommand{\CW}{{\cal W}}
\newcommand{\CX}{{\cal X}}
\newcommand{\CY}{{\cal Y}}
\newcommand{\CZ}{{\cal Z}}

\newcommand{\biU}{\hbox{$\bf\it U$}}
\newcommand{\rmn}[1]{{\rm (#1)\,\,}}

\newcommand{\bara}{\bar{a}}
\newcommand{\barA}{\bar{A}}
\newcommand{\barf}{\bar{f}}
\newcommand{\barh}{\bar{h}}
\newcommand{\bark}{\bar{k}}
\newcommand{\barn}{\bar{n}}
\newcommand{\barz}{\bar{z}}
\newcommand{\barkk}{\bar{\bar{k}}}
\newcommand{\barQ}{\bar{Q}}
\newcommand{\barga}{\bar{\ga}}
\newcommand{\bargep}{\bar{\gep}}
\newcommand{\bargb}{\bar{\gb}}
\newcommand{\bargl}{\bar{\gl}}
\newcommand{\bargd}{\bar{\gd}}
\newcommand{\bargo}{\bar{\go}}
\newcommand{\hatD}{\hat{D}}
\newcommand{\hatK}{\hat{K}}
\newcommand{\hatP}{\hat{P}}
\newcommand{\hatS}{\hat{S}}
\newcommand{\hatT}{\hat{T}}
\newcommand{\hatV}{\hat{V}}
\newcommand{\hatv}{\hat{v}}

\newcommand{\nuo}{\nu_{0}}

\newcommand{\veca}{\vec{a}}
\newcommand{\vecb}{\vec{b}}
\newcommand{\vecc}{\vec{c}}
\newcommand{\vecd}{\vec{d}}
\newcommand{\vecf}{\vec{f}}
\newcommand{\vecg}{\vec{g}}
\newcommand{\vecm}{\vec{m}}
\newcommand{\vecS}{\vec{S}}
\newcommand{\vecu}{\vec{u}}
\newcommand{\vecx}{\vec{x}}
\newcommand{\vecy}{\vec{y}}
\newcommand{\vecz}{\vec{z}}
\newcommand{\vecgf}{\vec{\gf}}
\newcommand{\vecgt}{\vec{\gt}}
\newcommand{\bgtu}{\bigtriangleup}
\newcommand{\tr}{\triangle}
\newcommand{\bgtd}{\bigtriangledown}
\newcommand{\bsl}{\backslash}
\newcommand{\nequ}{\equiv\!\!\!\!\! /}

\newcommand{\tl}{\tilde}
\newcommand{\tlb}{\tilde{b}}
\newcommand{\tlc}{\tilde{c}}
\newcommand{\tly}{\tilde{y}}
\newcommand{\tlG}{\tilde{G}}
\newcommand{\tlgt}{\tilde{\tau}}
\newcommand{\tlgw}{\tilde{\gw}}
\newcommand{\tlgs}{\tilde{\gs}}
\newcommand{\tlgm}{\tilde{\mu}}
\newcommand{\tlF}{\tilde{F}}
\newcommand{\tlJ}{\tilde{J}}
\newcommand{\tlgf}{\tilde{\phi}}
\newcommand{\tlY}{\tilde{Y}}
\newcommand{\tlgl}{\tilde{\lambda}}

\newcommand{\Raro}{\Rightarrow}
\newcommand{\LRaro}{\Leftrightarrow}
\newcommand{\raro}{\rightarrow}
\newcommand{\laro}{\leftarrow}
\newcommand{\Llraro}{\Longleftarrow}
\newcommand{\Lrraro}{\Longrightarrow}
\newcommand{\LLRraro}{\Longleftrightarrow}
\newcommand{\lngl}{\langle}
\newcommand{\rngl}{\rangle}

\newcommand{\bcap}{\bigcap}
\newcommand{\bcup}{\bigcup}
\newcommand{\sub}{\subset}
\newcommand{\cd}{\cdot}
\newcommand{\itms}[1]{\item[[#1]]}
\newcommand{\lnri}{\lim_{n\raro\infty}}
\renewcommand{\i}{\infty}

\newcommand{\ZZ}{{\mathchoice {\hbox{$\sf\textstyle Z\kern-0.4em Z$}}
{\hbox{$\sf\textstyle Z\kern-0.4em Z$}}
{\hbox{$\sf\scriptstyle Z\kern-0.3em Z$}}
{\hbox{$\sf\scriptscriptstyle Z\kern-0.2em Z$}}}}

\newcommand{\AAA}{{\mathchoice {\hbox{$\sf\textstyle A\kern-0.4em A$}}
{\hbox{$\sf\textstyle A\kern-0.4em A$}}
{\hbox{$\sf\scriptstyle A\kern-0.3em A$}}
{\hbox{$\sf\scriptscriptstyle A\kern-0.2em A$}}}}

\newcommand{\RR}{{\rm I\!R}}
\newcommand{\DD}{{\rm I\!D}}
\newcommand{\EE}{{\rm I\!E}}
\newcommand{\MC}{{\rm I\!\!\!C}}

\newcommand{\NN}{{\rm I\!N}}
\newcommand{\BB}{{\rm I\!B}}

\newcommand{\Cc}{{\mathchoice {\setbox0=\hbox{$\displaystyle\rm C$}\hbox
{\hbox
to0pt{\kern0.4\wd0\vrule height0.9\ht0\hss}\box0}}
{\setbox0=\hbox{$\textstyle\rm C$}\hbox{\hbox
to0pt{\kern0.4\wd0\vrule height0.9\ht0\hss}\box0}}
{\setbox0=\hbox{$\scriptstyle\rm C$}\hbox{\hbox
to0pt{\kern0.4\wd0\vrule height0.9\ht0\hss}\box0}}
{\setbox0=\hbox{$\scriptscriptstyle\rm C$}\hbox{\hbox
to0pt{\kern0.4\wd0\vrule height0.9\ht0\hss}\box0}}}}

\newcommand{\QQ}{{\mathchoice {\setbox0=\hbox{$\displaystyle\rm
Q$}\hbox{\raise
0.15\ht0\hbox to0pt{\kern0.4\wd0\vrule height0.8\ht0\hss}\box0}}
{\setbox0=\hbox{$\textstyle\rm Q$}\hbox{\raise
0.15\ht0\hbox to0pt{\kern0.4\wd0\vrule height0.8\ht0\hss}\box0}}
{\setbox0=\hbox{$\scriptstyle\rm Q$}\hbox{\raise
0.15\ht0\hbox to0pt{\kern0.4\wd0\vrule height0.7\ht0\hss}\box0}}
{\setbox0=\hbox{$\scriptscriptstyle\rm Q$}\hbox{\raise
0.15\ht0\hbox to0pt{\kern0.4\wd0\vrule height0.7\ht0\hss}\box0}}}}

\newcommand{\lun}[1]{\bigcup\limits_{#1}}

\newcommand{\II}{{\bf I\kern -1pt I}}

\newcommand{\half}{\frac{1}{2}}
\newcommand{\singcol}[2]{\left[\begin{rray}{c}#1\\#2\end{array}\right]}
\newcommand{\singcolb}[2]{\left(\begin{array}{c}#1\\#2\end{array}\right)}
\newcommand{\doubcol}[4]{\left[\begin{array}{cc}#1&#2\\#3&#4\end{array}\right]}
\newcommand{\doubcolb}[4]{\left(\begin{array}{cc}#1&#2\\#3&#4\end{array}\right)}
\newcommand{\sing}[1]{\langle #1\rangle}

\newcommand{\susim}{\stackrel{succ}{sim}}
\newcommand{\eqqdf}{\stackrel{def}{\equiv}}
\newcommand{\eqdf}{\stackrel{def}{=}}
\newcommand{\inj}{\stackrel{1-1}{\longrightarrow}}
\newcommand{\surj}{\vbox{\hbox{$\longrightarrow $
                  \kern -22pt \hbox{\lower 2.5pt  \hbox{\tiny onto}}
                  \kern -16pt \hbox{\raise 5pt  \hbox{\tiny 1-1}}
                  \kern 3pt}}}
\newcommand{\uarrow}[2]{\vbox{\hbox{$\longrightarrow $
                  \kern -16pt \hbox{\raise 5pt  \hbox{\tiny $#1$}}
                  \kern 10pt }}}
\newcommand{\AEQ}{$\forall \exists \:$ \,}
\newcommand{\A}{\forall}
\newcommand{\E}{\exists}

	\newcommand{\AH}{A\# H}

\newcommand{\R}{\hbox{I\kern-.1500em \hbox{\sf R}}}
\newcommand{\Q}
   {\hbox{${\rm Q} \kern -7.5pt \raise 2pt \hbox{\tiny$|$}\kern 7.5pt$}}
\newcommand{\C}
   {\hbox{${\rm C} \kern -7.5pt \raise 2pt \hbox{\tiny$|$}\kern 7.5pt$}}

\newcommand{\qd}{\kern 5pt\vrule height8pt width0.5pt depth0pt}

\newcommand{\bd}{\begin{description}}
\newcommand{\ed}{\end{description}}

\newcommand{\bc}{\begin{center}}
\newcommand{\ec}{\end{center}}

\newcommand{\beq}{\begin{equation}}
\newcommand{\eeq}{\end{equation}}

\newcommand{\ben}{\begin{enumerate}}
\newcommand{\een}{\end{enumerate}}

\newcommand{\seqn}[2]{\langle#1,\ldots ,#2\rangle}
\newcommand{\sseqn}[2]{\langle#1\,|#2\rangle}
\newcommand{\ssseqn}[2]{\langle#1\,|#2\rangle}
\newcommand{\setn}[2]{$\{\,$#1$\:|\:$#2$\}$}
\newcommand{\setm}[2]{\{\ #1\:|\:#2\}}

\newcommand{\setp}[2]{\{\,#1\: : \:#2\}}

\newcommand{\fsetn}[2]{\{\,#1,\ldots ,#2\}}
\newcommand{\ssetn}[2]{\Llbc\,#1\,|#2\Lrbc}
\newcommand{\sssetn}[2]{\LLlbc\,#1\,|#2\LLrbc}
\newcommand{\pair}[2]{\langle#1 ,#2\rangle}
\newcommand{\trpl}[3]{\langle#1 ,#2 ,#3 \rangle}
\newcommand{\fnn}[3]{#1:#2 \raro #3}
\newcommand{\surjfn}[3]{#1:#2 \surj #3}
\newcommand{\pgn}[3]{\langle#1 ,#2 ; \, #3 \rangle}
\newcommand{\iso}[3]{\m{ #1 : #2 \cong #3 }}

\newcommand{\hsetn}[2]{\hbox{\{\,#1\,|#2\}}}

\newcommand{\azc}{\hbox{$\aleph_{0}$-categorical }}
\newcommand{\azcy}{\hbox{$\aleph_{0}$-categoricity }}
\newcommand{\meq}{\hbox{$M^{Eq}$}}
\newcommand{\feq}{\hbox{$f^{Eq}$}}
\newcommand{\ut}{U-tree }
\newcommand{\uts}{U-trees }



\def\brd{\vrule height 2pt width .5pt depth 1.5pt}
\newcommand{\ol}{\overline}
\newcommand{\sm}{\setminus}
\newcommand{\minus}{^{-1}}
\newcommand{\bbraces}[1]{\left\{ #1 \right\}}


\def\newtheorems{\newtheorem{theorem}{Theorem}[section]
		 \newtheorem{cor}[theorem]{Corollary}
		 \newtheorem{prop}[theorem]{Proposition}
		 \newtheorem{lemma}[theorem]{Lemma}
		 \newtheorem{defn}[theorem]{Definition}
		 \newtheorem{claim}[theorem]{Claim}
		 \newtheorem{sublemma}[theorem]{Sublemma}
		 \newtheorem{example}[theorem]{Example}
		 \newtheorem{remark}[theorem]{Remark}
		 \newtheorem{question}[theorem]{Question}
		 \newtheorem{conjecture}{Conjecture}[section]}
\def\rar{\mathop{\longrightarrow}\limits}
\def\max{\mathop{\rm max}\limits}

\newcommand{\diaga}{\begin{picture}(40,12)
\put(24,20.2){\line(2,-1){69}}%
                   \end{picture}}

\newcommand{\thediaga}{\raisebox{-.2ex}{\makebox[0pt]{\diaga}}}

\newcommand{\rboxa}[1]{\raisebox{1.5ex}[2ex]{\hspace{.12in} #1}}

\newtheorem{theorem}{Theorem}
\newcounter{examp}
\font\bigcaps=cmcsc10 scaled 1250 %
\newcommand{\ExaLabel}[1]{\label{eqn:#1}}
\newcommand{\ExaRef}[1]{\ref{eqn:#1}}
\newcommand{\eqnref}[1]{(\ref{#1})}
\newcommand{\eop}[0]{\begin{flushright} $\blacksquare$ \end{flushright}}
\newtheorem{proposition}{Proposition}
\newtheorem{corollary}{Corollary}
\newtheorem{lemma2}{Lemma}
\newtheorem{remark}{Remark}
\newtheorem{example}{Example}
\newtheorem{question}{Question}
\newtheorem{definition}{Definition}
\newtheorem{property}{Property}
\newtheorem{conjecture}{Conjecture}

\newcommand{\norm}[1]{\left\lVert#1\right\rVert}
\newcommand{\BlackBox}{
\rule{1.1ex}{1.2ex} }
\newenvironment{Proof}[1]{\ifx\ShowProofs\True \addvspace{1.0ex}
{\bf Proof:}#1~\BlackBox\fi}{}
\def\proofi{\futurelet\next\lookforbracket}
\def\lookforbracket{\ifx\next[\let\go\usespecialterm
\else\let\go\relax
\vskip1sp\noindent{\bf Proof:}\enskip\relax\ignorespaces\fi\go}
\def\usespecialterm[#1]{\vskip12pt
\noindent{\bf #1:}\enskip\relax\ignorespaces}

\baselineskip 25pt
\title{\fontfamily{cmss} \selectfont Iterative process of order $2$ without inverting the derivative}
\vskip 15pt

\author{\setcounter{footnote}{6}Tamara Kogan \thanks{Amit Educational Network, Beer-Sheva, Israel, e-mail: tkogan37@gmail.com}$\;\;,$
\setcounter{footnote}{5}Ariel Sapir \thanks{Departments of Mathematics and Computer Science, Ben-Gurion University, Beer-Sheva, Israel}$\;\;,$
\setcounter{footnote}{0}Luba Sapir \thanks{Departments of Mathematics and Computer Science, Ben-Gurion University, Beer-Sheva, Israel,
e-mail: lsapir@bgu.ac.il}$\;\;,$
\setcounter{footnote}{5}Eytan Sapir \thanks{Departments of Mathematics and Computer Science, Ben-Gurion University, Beer-Sheva, Israel}$\;\;,$
\setcounter{footnote}{4}Amir Sapir \thanks{Department of Computer Science, Sapir Academic College, Sha'ar HaNegev, Yehudah, Israel and The Center for Advanced Studies in Mathematics, Ben-Gurion University, Beer-Sheva, Israel, e-mail: amirsa@cs.bgu.ac.il}}
\maketitle

\fontfamily{cmss} \selectfont
\abstract{We prove the sufficient conditions for convergence of a certain iterative process of order $2$ for solving nonlinear functional equations, which does not require inverting the derivative. We translate and detail our results for a system of nonlinear equations, and apply it for some numerical example which illustrates our theorems.}

\large
\baselineskip 16pt
\setlength\parindent{24pt}

\section{\textnormal{\bf \fontfamily{cmss} \selectfont Introduction}}~\label{intro}
\vskip -24pt
Many problems in Applied Mathematics and Engineering are solved by solution of nonlinear equations or its systems.
Usually, finding the exact solution of those equations is a very hard task at best, or even impossible sometimes. In these cases, we usually seek a numerical approximation of the solution by iterative methods, and Newton's method is one of the most popular, due to its simplicity.

However, Newton's method requires inverting the derivative on each step. For example, for solving systems of many nonlinear equations, this could mean inverting large size matrices, involving a large amount of computations.

Due to that, in 1974, one of the authors (Kogan \cite{Ko1974}) suggested an iterative method for a given nonlinear functional equation, which does not require inverting the derivative, starting from the second step. Yet, the order of this method is the same as that of Newton's, as will be proved in the sequel.

In this paper we also focus on applications of our results to the problem of solving systems of nonlinear equations. We provide two theorems, regarding sufficient conditions for convergence of our process and for existence of a solution to the system.
In the numerical example, our theoretical results are illustrated by computation and drawn in Figure (\ref{theorem1}).


\newpage

\subsection{\textnormal{\bf \fontfamily{cmss} \selectfont General Setup}}~\label{generalsetup}
We consider equations of the form:
\begin{equation}\label{eq1} 
P(x) = 0,
\end{equation}
where $P$ is a twice-derivable (in the sense of Frechet) operator, which transforms a Banach space $X$ to a space of the same type.

In {\cite{Ko1974}} was suggested the following iterative method:
\begin{equation}\label{iter1new} 
x_{n+1} = x_{n} - U_{n}P(x_{n}), \qquad\qquad n = 1,2,\ldots,
\end{equation}
where
$$U_{0} = [P^{\prime}(x_{0})]^{-1}$$
$$U_{n} = [2I - U_{n-1} P^{\prime}(x_{n})] U_{n-1} \qquad\qquad n=1,2,\ldots $$
and $I$ is the identity operator.

%
%
%
%

\subsection{\textnormal{\bf \fontfamily{cmss} \selectfont Structure of the Paper}}~\label{structure}
\qquad The paper is organized as follows: Section~\ref{res} describes several previous results.  Subsection~\ref{partA} presents a brief fundamental research, starting from  Kantorovich \cite{Kan1948a, Kan1948b} and up until modern time results.  Subsection~\ref{partB} provides one of the formulation of  Newton-Kantorovich Theorem.

Section~\ref{def} provides our main results and consists of two parts. In Subsection~\ref{main31} we state and provide the principal\footnote{A fully detailed proof can be found in the appendix. } ideas of the proof of Theorem \ref{theorem1}, which determines the order of convergence of this process and sufficient conditions for the existence of a solution in a certain domain of Banach space. The proof is based on the methodology of Kantorovich. In  Subsection~\ref{main32} we apply the suggested method for solving a system of nonlinear equations
\begin{equation}\label{system} 
f_{i}\left(x_1, \ldots, x_{n}\right) = 0, \qquad\qquad i=1,\ldots,n
\end{equation}
%

In Section \ref{num} we supply a numerical example. In this section we illustrate our main results (Theorem \ref{theorem2} and Theorem \ref{theorem3}) for proving
 the existence of solution of a certain system of equations and the convergence of our method to the solution.

%

Section \ref{conc} concludes the paper with a brief summary.

A full proof of our general theorem -- Theorem \ref{theorem1} -- is included in the appendix.

\newpage

\section{\textnormal{\bf \fontfamily{cmss} Previous Results}}~\label{res}
$\hskip 5pt$

\subsection{\textnormal{{\bf\fontfamily{cmss} \selectfont Iterative Methods for Functional Equations}}}~\label{partA}

One of the first fundamental researches of iterative processes appeared in 1948. In that research, L.V. Kantorovich \cite{Kan1948a,Kan1948b} obtained the range and the rate of convergence of Newton's method for nonlinear functional equation in Banach spaces and sufficient conditions for the existence of the solution.

By choosing specific spaces, Kantorovich managed to transform his results to solving of non-linear equations, system of equations and also differential and integral equations.

Following the methodology of L.V. Kantorovich, beginning in the 60's, many analogical results for other iterative processes have been found out (see, for example, \cite{Ser1961, Kogan1978}).

In recent decades, there was a significant advancement in tackling this problem, as can be seen in 
\cite{AmBu2007, AmBuPl2004, CoTo2011, JaSe2018, Kell2018, SoKaHa2018, Sidi2008, Anastasia2020, KoganSapir2007}.
An interesting book \cite{BuAm2016} by Sergio Amat and Sonia Busquier presents many results and advancements.

\subsection{\textnormal{{\bf\fontfamily{cmss} \selectfont Newton-Kantorovich Theoerem}}}~\label{partB}
\qquad Let the following conditions hold:\\
\begin{enumerate}
\item For an initial approximation $x_{0}$, there exists an operator $\Gamma_{0} = \left[P^{'}(x_{0})\right]^{-1}$  and $\norm{\Gamma_{0}}\leq B_0$.
\item $\norm{\Gamma_{0}\cdot P(x_{0})} \leq \eta_{0}$
\item $\norm{P''(x)}\leq K$ in the region
\begin{equation}\label{9obl}
 \norm{ x- x_0} \leq N\left( h_0 \right) \cdot \eta_0 = \frac{1-\sqrt{1-2h_0}}{h_0} \cdot \eta_0
\end{equation}
\item $h _0 = B_0 \cdot  \eta_0 \cdot K \le \frac{1}{2}.$
\end{enumerate}
Then equation~(\ref{eq1}) has the solution $x^{*}$ in the region~(\ref{9obl}) and
$$\norm{x^{(n)}- x^{*}}\leq \frac{1}{2^{n-1}}\cdot \left(2h_0\right )^{2^n-1} \eta_0.$$

\subsection{\textnormal{{\bf\fontfamily{cmss} \selectfont Non-linear systems of equations}}}~\label{systems}
\qquad We consider the following system of nonlinear equations:
\begin{equation}\label{non_linear_system1}
\left\{
\begin{array}{lll}
          f_1(x_1,x_2,\ldots x_n) &=& 0 \\
          f_2(x_1,x_2,\ldots x_n) &=& 0 \\
         \ldots &&  \\
           f_n(x_1,x_2,\ldots x_n) &=& 0,
        \end{array}
\right.\end{equation}

 which can be equivalently written in the form:
 \begin{equation}\label{sysP}
   P(\vec{X})=\vec{0},
  \end{equation}
where:
$P = \left( \begin{array}{c} f_1 \cr \ldots \cr \ldots \cr f_n
 \end{array} \right),$ $\vec{X} = \left( \begin{array}{c} x_1 \cr \ldots \cr \ldots \cr x_n
 \end{array} \right)$ and $\vec{0} = \left( \begin{array}{c} 0 \cr \ldots \cr \ldots \cr 0
 \end{array} \right)$.

The operator $P$ maps the $n$-dimensional space $X$ to itself. It is known that $P^{'}(\vec{X})$ is the Jacobian matrix. That is, $P^{'}(\vec{X})=\left(\frac{\partial f_i}{\partial x_j}\right)$
and $P^{''}(\vec{X})=\left(\frac{\partial^2 f_i}{\partial x_j \partial x_k}\right),$ $\;\; i,j,k=1, \ldots, n$ (See, for example, \cite{Kan1948a, Kan1948b}).

By choosing different metrics in the $n$-dimensional space, Kantorovich \cite{Kan1948a,Kan1948b} obtained different formulations of Newton-Kantorovich Theorem. In the sequel, we use his estimates (\ref{200}) -- (\ref{gamma0lambda}).\\
In the case $\norm{\vec{X}}=\max_{i} |x_i|$, he got:
 \begin{equation}\label{200}
   \norm{\Gamma_{0}}= \norm{\left(P^{'}(\vec{X}^{0})\right)^{-1}}\leq \max_{i} \frac{1}{|\Delta|} \sum\limits_{j=1}^n |A_{ij}|,
  \end{equation}
  and
  \begin{equation}\label{201}
   \norm{P^{''}(\vec{X})}\leq L\cdot n^2,
   \end{equation}
   where:
 \begin{itemize}
   \item[]$\Delta=\det\left(\frac{\partial f_i}{\partial x_j}\right)|_{\vec{X}=\vec{X}^{(0)}} $,
   \item[]$A_{ij}$ is $(i,j)$-cofactor  of $\left(\frac{\partial f_i}{\partial x_j}\right)|_{\vec{X}=\vec{X}^{(0)}} $,
   \item[]and $|\frac{\partial^2 f_i}{\partial x_j \partial x_k}| \leq L$ for any $1\leq i,j,k \leq n$ in the considered area.
 \end{itemize}
For   $\norm{\vec{X}}=\sqrt{\sum\limits_{i=1}^n x^2_i}$, he got:
\begin{equation}\label{202}
   \norm{\Gamma_{0}}\leq \left(\frac{1}{|\Delta|} \sum\limits_{i,k=1}^n A^2_{ik} \right)^\frac{1}{2},
  \end{equation}
and
\begin{equation}\label{203}
   \norm{P^{''}(\vec{X})}\leq L\cdot n \cdot \sqrt{n},
   \end{equation}
where $|\frac{\partial^2 f_i}{\partial x_j \partial x_k}| \leq L$ for any $1\leq i,j,k \leq n$ in the considered domain.

Moreover, if the eigenvalue $\lambda$ of a maximal modulus of the matrix $\Gamma_{0}\cdot \Gamma^T_{0}$ is known, then a more accurate estimate is:
\begin{equation}\label{gamma0lambda}
   \norm{\Gamma_0}\leq \sqrt{\lambda}
   \end{equation}

\newpage

 \section{\textnormal{\bf \fontfamily{cmss} \selectfont Main Results}}~\label{def} $\hskip 5pt$
\vskip -24pt

\subsection{\textnormal{{\bf\fontfamily{cmss} \selectfont Convergence of Iterative Process}}}~\label{main31}
Recall, that we solve Equation
\begin{equation}\label{eq1again} 
P(x) = 0,   \tag{\ref{eq1}}
\end{equation}
where $P(x)$ is a twice-derivable (in the sense of Frechet) operator which transforms a Banach space $X$ to a space of the same type. Assume that $U_{0} = [P^{\prime}(x_{0})]^{-1}$ exists, where $x_{0} \in X$ is an initial approximation to the solution of \eqnref{eq1}.

For solving \eqnref{eq1} we use the following iterative process:
\begin{equation}\label{iter1} 
x_{n+1} = x_{n} - U_{n}P(x_{n}), \qquad\qquad (n=0,1,\ldots), \tag{\ref{iter1new}}
\end{equation}
where:
\begin{equation}\label{detail1} 
U_{n} = [2I - U_{n-1} P^{\prime}(x_{n})] U_{n-1},
\end{equation}
and $I$ is the identity operator.

We show that Process (\ref{iter1}) converges at the same order as Newton-Kantorovich, yet requires only a single inversion of $P^{\prime}(x_{0})$.

\subsubsection{\textnormal{{\large\bf\fontfamily{cmss} \selectfont Sufficient Conditions for Convergence}}}~\label{main311}
\begin{theorem}~\label{theorem1}
 Let $a$ be a real root of the equation:
$$a^3 + 2a^2 +3a -2 = 0, \qquad 0.477 < a < 0.478$$.
Assume that for $B, \eta, K$ it holds that:
\begin{itemize}
\item[$1^{\circ}$] $\norm{U_{0}} \le B; \,\, \norm{P(x_{0})} \le \eta$.
\item[$2^{\circ}$] $\norm{P^{\prime\prime}(x)} \le K$, for any $x\in G$, where
\begin{equation}\label{denoteSN}
\begin{aligned}
G = \left\{ x \left| \norm{x - x_{0}} \le \frac{2-a-a^2}{2(1-a-a^2)}B\eta  \right.\right\}
\end{aligned}
\end{equation}
\item[$3^{\circ}$] $h = B^{2} \eta K \le a.$
\end{itemize}

Then \eqnref{eq1} has a solution $x^{*}$ in the region $G$, the sequence of \eqnref{iter1} converges and:
\begin{equation}\label{inTheorem} 
\norm{x_{n} - x^{*}} = \frac{h(1+h)}{2} \frac{S^{n-1}}{1 - S} (N_{1}h)^{2^{n}-2}B\eta
\end{equation}
\end{theorem}
where:
\begin{equation*}\label{denoteSNonly}
\begin{aligned}
S = \frac{2(1+h)}{h^2 + 2h + 3}; \,\, N_{1} = \frac{h^2 + 2h + 3}{2}
\end{aligned}
\end{equation*}

\subsubsection{\textnormal{{\large\bf\fontfamily{cmss} \selectfont Sketch of Convergence Proof}}}~\label{main312}
In this subsection we will show only the principal ideas of the proof. This is intended for readers who are interested in the practical aspects of the question. The detailed proof appears in the appendix.

We show that, for any $k = 1,2,3 ,\ldots$, the following inequalities:
\begin{equation}\label{norm x_k} 
\norm{x_{k} - x_{k-1}} \le \alpha_{k}h^{2^{k-1}-1}B\eta
\end{equation}

\begin{equation}\label{norm P x_k} 
\norm{P(x_{k})} \le \beta_{k}h^{2^{k}-1}\eta
\end{equation}

\begin{equation}\label{diff J P x_k} 
\norm{I - P^{\prime}(x_{k})U_{k-1}} \le A_{k}h^{2^{k-1}}
\end{equation}

\begin{equation}\label{U_k} 
\norm{U_{k}} \le c_{k}B
\end{equation}

\begin{equation}\label{norm x_k to x_0} 
\norm{x_{k} - x_{0}} \le \gamma_{k}B\eta
\end{equation}
hold, where $\alpha_{k}, \beta_{k}, A_{k}, c_{k}, \gamma_{k}$ are bounded constants.

From these inequalities we obtain by induction, that all the approximations $x_{k} \in G$ and $\{x_{k}\}$ is a fundamental sequence.

By properties of Banach Space and fundamentality of $\{x_{k}\}$ emanates existence of  $x^{*} = \lim\limits_{n \rightarrow \infty} x_{n}$, which is the solution of $(\ref{eq1})$.

\subsection{\textnormal{{\Large \bf\fontfamily{cmss} \selectfont Nonlinear systems of equations}}}~\label{main32}
\qquad Here we apply our method to solve the following system of nonlinear equations:
\begin{equation}\label{non_linear_equations}
\left\{
\begin{array}{lll}
          f_1(x_1,x_2,\ldots x_n) &=& 0 \\
          f_2(x_1,x_2,\ldots x_n) &=& 0 \\
         \ldots &&  \\
           f_n(x_1,x_2,\ldots x_n) &=& 0,
        \end{array}
\right. \tag{\ref{non_linear_system1}}
\end{equation}

As we mention in~\ref{systems} it can be equivalently written in the form:
\begin{equation}
   P(\vec{X})=\vec{0}\tag{\ref{sysP}},
\end{equation}

where
$P = \left( \begin{array}{c} f_1 \cr \ldots \cr \ldots \cr f_n
 \end{array} \right),$ $\vec{X} = \left( \begin{array}{c} x_1 \cr \ldots \cr \ldots \cr x_n
 \end{array} \right)$ and $\vec{0} = \left( \begin{array}{c} 0 \cr \ldots \cr \ldots \cr 0
 \end{array} \right)$.

 The suggested iterative process can be written as following:
 \begin{equation}\label{ourpross}
   \vec{X}^{(k+1)}=\vec{X}^{(k+1)}-U_k P(\vec{X}^{(k)}), \qquad k=0,1,2,\cdot \;.
  \end{equation}

The matrix $U_{k+1}$ can be expressed by the recursive expression $$U_{k+1}=\left[2I-U_k P^{'}(\vec{X}^{(k+1)})\right]U_k,$$
where: $\vec{X}^{(0)}$  is an initial approximation to the solution $\vec{X}^{*}$, $P^{'}(\vec{X})=\left(\frac{\partial f_i}{\partial x_j}\right)$,\\ $U_0=\left(P^{'}(\vec{X}^{(0)})\right)^{-1}$
 and $I$ is the identity matrix.

\begin{equation*}
P^{''} \left(\vec{X} \right) = \left( \frac{\partial^2 f_i}{\partial x_j \cdot \partial x_k} \right), \qquad 1\leq i,j,k \leq n
\end{equation*}

The operator $P$ maps the $n$-dimensional space $X$ to itself.


Choosing distinct metrics for the $n$-dimensional space, we obtain 2 theorems.

The first case is $\norm{\vec{X}}=\max_{i} |x_i|$. In this case, we obtain:

 \begin{theorem}~\label{theorem2}
If the following conditions hold for $B, \eta, K$:
\begin{itemize}
\item[$1^{\circ}$] $ \max_{i} \left| f_i \left( \overrightarrow{X}^{\left( 0 \right)} \right) \right| \leq \eta $,
\item[$2^{\circ}$] $\frac{1}{\left| \Delta \right|} \cdot \sum\limits_{i,k=1}^{n} \left| A_{i,k} \right| \leq B$ (recall that $A_{i,k}$ is $(i,k)$-cofactor  of $\left(\frac{\partial f_i}{\partial x_k}\right)|_{\vec{X}=\vec{X}^{(0)}}$),
\item[$3^{\circ}$] $n^{2} \cdot L \leq K $ for $L=\max_{i,j,k}{ \left| \frac{\partial^2 f_i}{\partial {x_j}\cdot \partial {x_k}} \right| }$ in our region,
\item[$4^{\circ}$] $B^2 \cdot \eta \cdot K \leq a $ (where $ 0.477 < a < 0.478$).
\end{itemize}

Then the iterative process \eqnref{iter1new} converges to the solution $\overrightarrow{X} ^{*}$ of the system, and
$$\norm{\overrightarrow{X}^{\left( n \right)} - \overrightarrow{X}^{*}} = \max_i \left| \overrightarrow{x_i}^{\left( n \right)} - \overrightarrow{x}_{i}^{*} \right| \leq \frac{h(1+h)}{2} \frac{S^{n-1}}{1 - S} (N_{1}h)^{2^{n}-2}B\eta.$$
\end{theorem}

The second case we consider is when  $\norm{\vec{X}}=\sqrt{\sum\limits_{i=1}^n x^2_i}$. Here, we obtain:
\begin{theorem}~\label{theorem3}
If the following conditions hold for  $B, \eta, K$:
\begin{itemize}
\item[$1^{\circ}$] $ \sqrt{\sum\limits_{i=1}^n x^2_i} \leq \eta $,
\item[$2^{\circ}$] $\left| U_0\right|\leq \sqrt{\lambda} = B$, where $\lambda$ is the eigenvalue\footnote{If finding $\lambda$ is a complicated task, we can require $\norm{U_0} \leq \sqrt{\left( \frac{1}{\left|\Delta\right|} \cdot \sum\limits_{i,k} {A_{i,k}}^2 \right)} $ instead.}  of a maximal modulus of the matrix $U_{0}\cdot U^T_{0}$.
\item[$3^{\circ}$] $\norm{P^{''} \left(x \right)} \leq  n \cdot \sqrt{n} \cdot L$,
\item[$4^{\circ}$] $B^2 \cdot \eta \cdot K \leq a $ (where $ 0.477 < a < 0.478$).
\end{itemize}

Then the iterative process \eqnref{iter1new} converges to the solution $\overrightarrow{X} ^{*}$ of the system, and
$$\norm{\overrightarrow{X}^{\left( n \right)} - \overrightarrow{X}^{*}} = \sqrt{ \sum\limits_{i=1} ^{n} \left( \overrightarrow{x_i}^{\left( n \right)} - \overrightarrow{x}_{i}^{*} \right)^{2} }\leq \frac{h(1+h)}{2} \frac{S^{n-1}}{1 - S} (N_{1}h)^{2^{n}-2}B\eta.$$
\end{theorem}
For the above theorems we employed inequalities \eqnref{200}, \eqnref{201}, \eqnref{202} and \eqnref{203}.

\newpage

\section{\bf\fontfamily{cmss} \selectfont Numerical Results}~\label{num}
The following example illustrates the suggested method.

Consider the system of non-linear equations:
\begin{equation}\label{nm_res_equations}
\left\{
\begin{array}{ll}
f(x_1,x_2)=2{x_1}^3 - {x_2}^2 - 1 = 0 \\
f(x_1,x_2)=x_1\cdot {x_2}^3 - x_2 - 4 =0
        \end{array}
\right.
\end{equation}
whereas the solution is in the domain:
\begin{equation}\label{nm_res_domain}
\left\{\begin{array}{ll}
0\leq x_1\leq 1.3 \\
0\leq x_2\leq 1.8
        \end{array}
\right\}
\end{equation}

In the case where the domain is unknown, one can construct it by using the algorithm in \cite{Kogan1964}, or by any other mean.

Let $\overrightarrow{X}^{(0)} = \left( \begin{matrix}
1.2\\
1.7\\
\end{matrix} \right)$.

Let us now verify the conditions of Theorem \ref{theorem2}.

\begin{enumerate}
  \item[1$^{\circ}$] $f_1\left( \overrightarrow{X}^{(0)} \right) = - 0.434, \,\, f_2\left( \overrightarrow{X}^{(0)} \right) = 0.1956$ \\
      $\eta = \max\left\{ 0.434, 0.1956 \right\} = 0.434$ \\
  \item[2$^{\circ}$] $P^{'} \left( \overrightarrow{X}\right) = \left( \begin{matrix}  6{x_1} ^2 & -2x_2 \\{x_2}^3 & 3x_1 {x_2}^2 -1 \\ \end{matrix} \right) $, $\,\, P^{'} \left( \overrightarrow{X} ^{(0)}\right) = \left( \begin{matrix}  8.64 & -3.4 \\4.913 & 9.404 \\ \end{matrix} \right) $ \\ \\
      $\Delta = 97.95 \neq 0 $, $\, B=\frac{9.404+4.913+3.4+8.64}{97.95}\leq 0.27$ \\
  \item[3$^{\circ}$] $\frac{\partial^2 f_1}{\partial {x_1}^2} = 12x_1$, $\frac{\partial^2 f_1}{\partial x_1 \partial x_2} = 0$, $\frac{\partial^2 f_1}{\partial {x_2}^2} = -2$ \\ \\
      $\frac{\partial^2 f_2}{\partial {x_1}^2} = 0$, $\frac{\partial^2 f_2}{\partial x_1 \partial x_2} = 3{x_2}^2 $, $\frac{\partial^2 f_2}{\partial {x_2}^2} = 6x_1 x_2$ \\ \\
      These functions are increasing in the domain $D$. Thus, we obtain        $$L=\max_{i,j,k=1,2} {\left| \frac{\partial^2 f_i}{\partial x_j \cdot \partial x_k} \right|}=\left|\frac{\partial^2 f_1}{\partial {x_1}^2}\right|=12\cdot 1.3 = 15.6$$
  \item[4$^{\circ}$] $h=B^2 \cdot \eta \cdot L\cdot n^2 = 0.27^2 \cdot 0.434\cdot 15.6\cdot 4 > a$
\end{enumerate}

This means that the conditions of Theorem \ref{theorem2} do not hold. Hence, it cannot be applied to determine whether a solution exists or the process converges.

We turn to check whether we can apply  Theorem \ref{theorem3}. Let us now verify its conditions:

\begin{enumerate}
  \item $\eta = \sqrt{0.434^2 + 0.1956^2}=0.476$ \\
  \item Since $\norm{U_0} = \sqrt{\lambda} $, where $\lambda$ is an eigenvalue of the matrix $U_0 \cdot {U_0}^{T}$: \\ \\
      $ U_0 \cdot {U_0} ^{T} = \begin{pmatrix} 0.010421 & -0.001753\\ -0.001753 & 0.010295\protect\end{pmatrix} $\\
      From the quadratic equation
      $$\left( 0.010421 - \lambda \right) \cdot \left( 0.010295 - \lambda \right) - 0.001753^2 = 0 $$
      we obtain the following eigenvalues:
      $$\lambda_1 = 0.0121, \lambda_2 = 0.0086,$$

      $$and \quad B=\sqrt{0.0121} = 0.11.$$
   \item $ K \leq n\sqrt{n} \cdot L = 2\sqrt{2} \cdot 15.6 < 44.1235 $.
   \item $ h = B^2 \cdot K \cdot \eta < 0.11^2 \cdot 44.1235 \cdot 0.476 < a $. \\
\end{enumerate}

This means that the conditions for Theorem \ref{theorem3} hold.

Thus, the equation has a solution $\overrightarrow{X} ^{*} = \begin{pmatrix} {x_1}^*\\ {x_2}^*\end{pmatrix} $ in the circle $G_0$ centered around $\begin{pmatrix} {1.2}\\ {1.7}\end{pmatrix}$  with the radius $r_0 = \frac{2-a-a^2}{2\cdot \left( 1-a-a^2 \right)}\cdot B \cdot \eta = 0.115$:
$$G_0 = \left\{ \begin{pmatrix} {x_1}\\ {x_2}\end{pmatrix} \,\middle| \, \left(x_1 -1.2\right)^2 + \left(x_2-1.7\right)^2 \leq {r_0}^2 \right\}$$

Yet, there exists a small probability of $\overrightarrow{X} ^{*}$ not being inside $D$ (see Figure \ref{theorem1}).

To overcome this, we tighten our region by computing $\overrightarrow{X} ^{(1)}$ and its surrounding circle $G_1$:\\
$ \overrightarrow{X} ^{(1)} = \begin{pmatrix} {x_1}\\ {x_2}\end{pmatrix} = \begin{pmatrix} 1.23488\\ 1.660982\end{pmatrix}$, \\
$ G_1 = \left\{ \begin{pmatrix} {x_1}\\ {x_2}\end{pmatrix} \Bigg| \, \left(x_1 - 1.23488\right)^2 +\left( x_2 - 1.660982\right)^2 \leq 0.028^2 \right\}$.

It can be seen that $\overrightarrow{X} ^{*} \in G_1 \subset D$.

\begin{figure}[ht]~\label{fig1}
    \begin{center}
    \includegraphics[width=10.8cm]{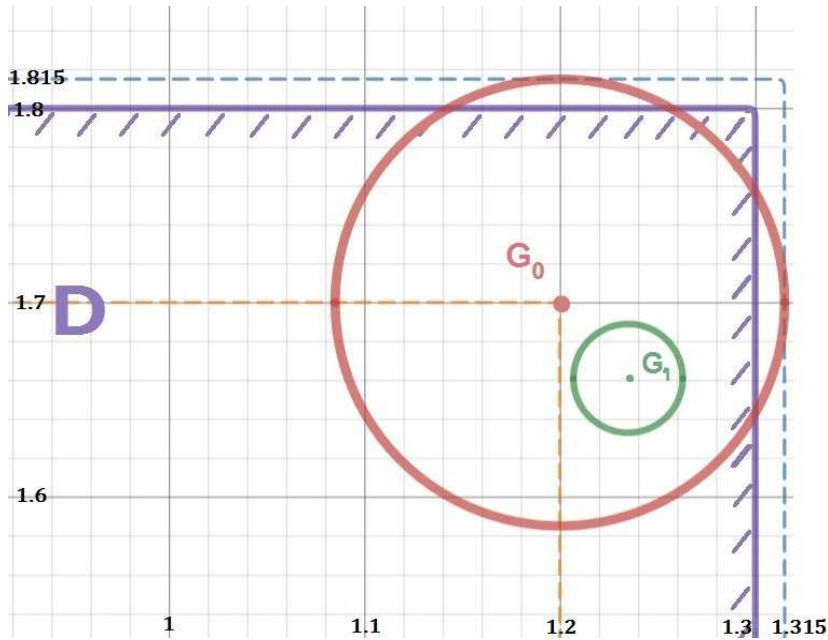}
    \caption{Illustration of $D$, $G_0$ and $G_1$.}
    \end{center}
\end{figure}
\newpage

\newcommand{\rbox}[1]{\raisebox{1ex}[2.3ex]{\hspace{.005in} #1}}

The values of $\overrightarrow{X} ^{(i)}$ and $P\left( \overrightarrow{X} ^{(i)} \right)$, for $i = 1,2,3,4$, are given in the table below:
\begin{table}[htbp]~\label{tb0}
\begin{adjustbox}{width=\columnwidth,center}
  \normalsize
  \begin{tabular}{|l|c|c|}
    \hline
    $i$ & $\overrightarrow{X} ^{(i)}$ &$P\left( \overrightarrow{X} ^{(i)} \right)$  \\
    \hline
    $0$ & $\begin{pmatrix} 1.2\\1.7\end{pmatrix}$ & $\begin{pmatrix} -4.34\cdot 10^{-1}\\1.1956\cdot 10^{-1}\end{pmatrix}$  \\
    \hline
    $1$ & $\begin{pmatrix} 1.234876263286\\1.660979680824\end{pmatrix}$ & $\begin{pmatrix} 7.32\cdot 10^{-3}\\2.28\cdot10^{-3}\end{pmatrix}$  \\
    \hline
   $2$ & $\begin{pmatrix} 1.234275470964\\1.661525517833\end{pmatrix}$ & $\begin{pmatrix} 1.217\cdot10^{-5}\\4.225\cdot 10^{-6}\end{pmatrix}$   \\
    \hline
    $3$ & $\begin{pmatrix} 1.234274484119\\1.661526466792\end{pmatrix}$ & $\begin{pmatrix} 5.287\cdot10^{-11}\\-1.875\cdot 10^{-11}\end{pmatrix}$   \\
    \hline
    $4$ & $\begin{pmatrix} 1.234274484114\\1.661526466796\end{pmatrix}$ & $\begin{pmatrix} -8.882\cdot10^{-20}\\0\end{pmatrix}$  \\
    \hline
  \end{tabular}
  \end{adjustbox}
\end{table}

\newpage

\section{\textnormal{\fontfamily{cmss} \selectfont Summary}}~\label{conc}
\vskip -24pt
We conclude with the following points:
\begin{enumerate}
\item The proposed iterative procedure has the same convergence order as Newton's Method, yet it differs in that, starting from Step 2, it does not require inverting the matrix. This reduces the amount of calculations.

\item For practical usage, one may often encounter a scenario in which checking the conditions for convergence is more of a complicated task than calculating the solution's approximations. In such cases, the authors recommend to compute $k$ approximations (for $k\geq 1$) and if $P\left(x_k \right)$ is small enough, one may check the conditions for existence of a solution with an initial approximation $x_k$. A small size ( ($\norm{P\left( x_k\right)} \leq  \eta$) ) of the parameter $\eta$  increases the probability of fulfilment of Theorem \ref{theorem2}'s or Theorem \ref{theorem3}'s  conditions.

\item In the case of a single equation  $f\left( x \right)=0$, the efficiency index of our procedure is equal to that of Newton's Method $I=\sqrt 2$. The efficiency index can be increased ($I=2$) if we build a nonstationary process of the same order of convergence as proposed in \cite{KoganSapir2020}.

 \end{enumerate}

\bibliographystyle{elsart-num-sort}
\baselineskip 10pt
\normalsize

\section{\textnormal{\bf \fontfamily{cmss} \selectfont Appendix}}~\label{appen}
\Large
\baselineskip 16pt
\setlength\parindent{24pt}
\begin{proofi}[Proof of Theorem \ref{theorem1}]
Define the numbers $\alpha_{k},\,\,\beta_{k},\,\,A_{k},\,\,c_{k}$ by the formulae
$$\alpha_{1} = 1, \qquad \beta_{1} = \frac{1}{2}, \qquad c_{1} = 1+h, \qquad A_{1} = 1$$

\begin{equation}\label{alpha_k} 
\alpha_{k} = c_{k-1}\beta_{k-1}
\end{equation}

\begin{equation}\label{beta_k} 
\beta_{k} = A_{k-1}^{2}\beta_{k-1} + \frac{1}{2}\alpha_{k}^{2}
\end{equation}

\begin{equation}\label{A_k} 
A_{k} = A_{k-1}^{2} + \alpha_{k}c_{k-1}
\end{equation}

\begin{equation}\label{c_k} 
c_{k} = c_{k-1}\left(1 + A_{k}h^{2^{k-1}}\right)
\end{equation}

Let
\begin{equation}\label{epsilon_k} 
\varepsilon_{k} = \alpha_{k}h^{2^{k-1}-1}, \,\, q_{k} = A_{k}h^{2^{k-1}}, \,\, \gamma_{k} = \sum\limits_{i=1}^{k}\varepsilon_{i}
\end{equation}

We will show by induction, that for any $k = 1,2,3 ,\ldots$, the following inequalities
\begin{equation}\label{norm x_k_fproof} 
\norm{x_{k} - x_{k-1}} \le \alpha_{k}h^{2^{k-1}-1}B\eta  \tag{\ref{norm x_k}}
\end{equation}

\begin{equation}\label{norm P x_k_fproof} 
\norm{P(x_{k})} \le \beta_{k}h^{2^{k}-1}\eta             \tag{\ref{norm P x_k}}
\end{equation}

\begin{equation}\label{diff J P x_k_fproof} 
\norm{I - P^{\prime}(x_{k})U_{k-1}} \le A_{k}h^{2^{k-1}} \tag{\ref{diff J P x_k}}
\end{equation}

\begin{equation}\label{U_k_fproof} 
\norm{U_{k}} \le c_{k}B                                  \tag{\ref{U_k}}
\end{equation}

\begin{equation}\label{norm x_k to x_0_fproof} 
\norm{x_{k} - x_{0}} \le \gamma_{k}B\eta                 \tag{\ref{norm x_k to x_0}}
\end{equation}
hold, and all the approximations $x_{k}$ are inside the region
$G$. \\

In the case of $k=1$, we obtain from \eqnref{iter1}, \eqnref{detail1} and $1^{\circ} - 3^{\circ}$ that:

$$\norm{x_{1} - x_{0}} \le \norm{U_{0}} \norm{P(x_{k})} \le B\eta = \alpha_{1}B\eta .$$

Since $1 \le \frac{2-a-a^2}{2(1-a-a^2)}$, we have that $x_{1}$ lies inside $G$ \\

$$ \norm{P(x_{1})} = \norm{P(x_{1}) - P(x_{0}) -  P^{\prime}(x_{0})(x_{1} - x_{0})} \le \frac{1}{2}h\eta = \beta_{1}h\eta $$

$$ \norm{I - P^{\prime}(x_{k})U_{0}} \le  \norm{P^{\prime}(x_{0}) - P^{\prime}(x_{1})} \cdot\norm{U_{0}} \le K B^{2} \eta = h = A_{1}h $$

$$ \norm{U_{1}} \le \norm{I + U_{0}\left(P^{\prime} (x_{0}) - P^{\prime} (x_{1})\right)} \cdot \norm{U_{0}} \le (1+h)B = c_{1}B $$
\\

Assume that for $k\leq n-1$ the conditions \eqnref{norm x_k} to \eqnref{norm x_k to x_0} hold and also that $\norm{x_{k} - x_{0}} \le \frac{2-a-a^2}{2(1-a-a^2)}B\eta .$

We will show that the same conditions apply for $k = n$. Indeed, from (\ref{iter1}), \eqnref{epsilon_k}, \eqnref{norm P x_k}, \eqnref{U_k}, \eqnref{norm x_k to x_0} we have:

$$ \norm{x_{n} - x_{n-1}} \le \norm{U_{n-1}} \cdot \norm{P(x_{n-1}) } \le c_{n-1}\beta_{n-1}h^{2^{n-1}-1}B\eta = \alpha_{n}h^{2^{n-1}-1}B\eta $$
$$ \norm{x_{n} - x_{0}} \le \norm{x_{n} - x_{n-1}} + \norm{x_{n-1} - x_{0}} \le \left(\gamma_{n-1} + \alpha_{n}h^{2^{n-1}-1}\right) B\eta = \gamma_{n}B\eta \, , $$
i.e. conditions \eqnref{norm x_k} and \eqnref{norm x_k to x_0} hold also for $k = n$.

We now turn to prove several inequalities that will be later used in order to prove \eqnref{norm P x_k} - \eqnref{U_k}.

Multiplying the two sides of \eqnref{A_k} by $h^{2^{k-1}}$, we obtain:
\begin{equation}\label{q_k_eq} 
 q_{k} = q_{k-1}^{2} + c_{k-1}\alpha_{k}h^{2^{k-1}}
\end{equation}

Employing (\ref{alpha_k}), (\ref{beta_k}), (\ref{c_k}) leads to:
$$ \alpha_{k}c_{k-1} = \left(1 + q_{k-1}\right)^{2} \alpha_{k-1}c_{k-2} \left[A_{k-2}^{2} + \frac{1}{2}\alpha_{k-1}c_{k-2} \right] $$

Multiplying the two sides of the above equality by $h^{2^{k-1}}$ and employing \eqnref{q_k_eq}, we obtain:
$$ \alpha_{k}c_{k-1}h^{2^{k-1}} = q_{k} - q_{k-1}^{2} \, , \qquad \alpha_{k-1}c_{k-2}h^{2^{k-2}} = q_{k-1} - q_{k-2}^{2} \, ,$$
, which leads to the relation:
$$ q_{k} = q_{k-1}^{2} + \frac{1}{2} \left( 1+q_{k-1} \right)^{2} \left( q_{k-1}^{2} - q_{k-2}^{4} \right) \, ,$$
, from which it emanates that:
\begin{equation}\label{q_k_leq} 
q_{k} \le N_{k-1}q_{k-1}^{2} \, ,
\end{equation}
, where $N_{k-1} = 1 + \frac{1}{2} \left( 1+q_{k-1} \right)^{2}.$


Recall that $q_{1} = h, \qquad N_{1}q_{1} = \frac{h^{3} + 2h^{2} + 3h}{2}, \qquad 0 \le h \le a$. Due to the increase of $\varphi(h) = h^{3} + 2h^{2} + 3h$ in $[0,a]$
from $a^{3} + 2a^{2} + 3a - 2=0$, it emanates that $N_{1}q_{1} \le 1$ always holds.

By induction, \eqnref{q_k_leq} employs the following property: The sequences $\{N_{n}\}, \{q_{n}\}$ are decreasing and $N_{k}q_{k} \le 1$ for any natural $k$.

Substituting $k=n$ and then $k=n-1$ in (\ref{q_k_eq}) together with condition (\ref{epsilon_k}) we get:
$$ q_{n} = q_{n-1}^{2} + c_{n-1}h\varepsilon_{n} \, , \qquad q_{n-1} = q_{n-2}^{2} + c_{n-2}h\varepsilon_{n-1}, $$
, leading to:
$$ \varepsilon_{n} = \varepsilon_{n-1} \frac{q_{n} - q_{n-1}^{2}}{\left(1+q_{n-1}\right)\left(q_{k-1} - q_{k-2}^{2}\right)} \, .$$

Since:
$$ \frac{q_{n} - q_{n-1}^{2}}{\left(q_{n-1} - q_{n-2}^{2}\right)} = \frac{1}{2} \left( 1+q_{n-1} \right)^{2} \left(q_{n-1} + q_{n-2}^{2}\right) \, ,$$
, we obtain:
$$ \varepsilon_{n} = \frac{1}{2}\varepsilon_{n-1}\left( 1+q_{n-1} \right) \left(q_{n-1} + q_{n-2}^{2}\right)$$

From the above equality and condition (\ref{q_k_leq}) we obtain:
$$ \varepsilon_{n} < \frac{1}{2}\left( 1+q_{n-1} \right) \left( 1+N_{n-1} \right) q_{n-2}^{2} \varepsilon_{n-1} < M q_{n-2}^{2} \varepsilon_{n-1}$$
, where $M = \frac{1}{2}\left( 1+h \right) \left( 1+N_{1} \right)$.

By \eqnref{q_k_leq} and the above inequality, we get that for $n=3,4,5 \ldots$ it holds that:
\begin{eqnarray*}
\begin{tabular}{rclcl}
$\varepsilon_{3}$& $<$& $M \varepsilon_{2} q_{1}^{2}$ &  $=$ & $\frac{M}{N_{1}^{2}} \varepsilon_{2} \left( N_{1}q_{1} \right)^{2}$\\
$\varepsilon_{4}$ &$<$ &$M^{2} \varepsilon_{2} N_{1}^{2}q_{1}^{6}$ &$=$ & $\left(\frac{M}{N_{1}^{2}}\right)^{2} \varepsilon_{2} \left( N_{1}q_{1} \right)^{6}$\\
$. $& & $.$ & &\\
$.  $& & $.$ & &\\
$.   $& & $.$ & &\\
\end{tabular}
\end{eqnarray*}

\begin{equation}\label{epsilon_n_leq} 
\begin{tabular}{rclcl}
$\varepsilon_{n}$& $<$ &$\left(\frac{M}{N_{1}^{2}}\right)^{n-2} \varepsilon_{2} \left(N_{1}q_{1}\right)^{2^{n-1}-2}$ & \qquad\qquad\qquad &
\end{tabular}
\end{equation}
$$\frac{M}{N_{1}^{2}} = \frac{\left( 1+q_{1} \right) \left( 1+N_{1} \right)}{2N_{1}^{2}} < \frac{1}{1+\frac{1+q_{1}^{2}}{2(1+q_{1})}} = S < 1 \, .$$

Note that $\varepsilon_{2} = \frac{1}{2} h(1+h)$ and $q_{1} = h$. Hence, we obtain:
$$ \varepsilon_{n} \le S^{n-2} \cdot \frac{h(1+h)}{2} \cdot \left(N_{1} h\right)^{2^{n-1}-2} $$

Denote $\gamma = \sum\limits_{n=1}^{\infty} \varepsilon_{n} = \lim\limits_{n \to \infty} \gamma_{n}$. Inequality \eqnref{epsilon_n_leq} emanates that:
$$\gamma = 1 + \varepsilon_{2}\sum\limits_{i=2}^{\infty} \frac{\varepsilon_{i}}{\varepsilon_{2}}
\le 1 + \varepsilon_{2}\sum\limits_{i=2}^{\infty} S^{i-2} \left(N_{1}h\right)^{2^{i-1}-2}
\le 1 + \frac{\varepsilon_{2}}{1 - S\left(N_{1}h\right)^{2}}
 $$
As we recall that:
$$ a^{3} + 2a^{2} + 3a - 2 = 0 \, ,$$
, we obtain:
$$ S\left(N_{1}h\right)^{2}
= \frac{h^{2}(h^{2} + 2h + 3)(1+h)}{2}
\le \frac{a(a^{3} + 2a^{2} + 3a)(1+a)}{2}
= a(1+a) \, ,$$
By the above inequality we get that:
$$ 1 \le \gamma \le \frac{2-a-a^2}{2(1-a-a^2)} $$
, and thus:
$$\norm{x_{n} - x_{0}} \le \gamma_{n}B\eta < \gamma B\eta \le \frac{2-a-a^2}{2(1-a-a^2)}B\eta$$

These inequalities will suffice for us. Let us resume to proving  \eqnref{norm P x_k} - \eqnref{U_k}, for 
$k=n$.

Since:
$$ I - P^{\prime}(x_{n-1})U_{n-1} =  I - P^{\prime}(x_{n-1})U_{n-2} - P^{\prime}(x_{n-1})\left(U_{n-1} - U_{n-2}\right) = [I - P^{\prime}(x_{n-1})U_{n-2}]^{2} \, ,$$
, we obtain that:
$$ \norm{I - P^{\prime}(x_{n-1})U_{n-1}} \le A_{n-1}^{2}h^{2^{n-1}} $$

Recalling that $x_{n} - x_{n-1} = -U_{n-1}P(x_{n-1})$, we obtain:
$$ \norm{P(x_{n-1}) + P^{\prime}(x_{n-1})\left(x_{n} - x_{n-1}\right)} \le \norm{I - P^{\prime}(x_{n-1})U_{n-1}} \cdot \norm{P(x_{n-1})} \le A_{n-1}^{2} \beta_{n-1} h^{2^{n}-1} \eta \, , $$

By using Taylor's formula, we then receive that:
\begin{eqnarray*}
\begin{tabular}{rcl}
$\norm{P(x_{n-1})}$ &$\le$ & $\norm{P(x_{n-1}) + P^{\prime}(x_{n-1})\left(x_{n} - x_{n-1}\right)} + \frac{1}{2}K\norm{x_{n} - x_{n-1}}^{2}$\\
\\
 &$\le $ & $\left(A_{n-1}^{2} \beta_{n-1} + \frac{1}{2}\alpha_{n}^{2}\right)
h^{2^{n}-1}\cdot  \eta$\\
\\
 & $ = $ & $ \beta_{n} h^{2^{n}-1} \eta$
\end{tabular}
\end{eqnarray*}

By the inequality:
$$I - P^{\prime}(x_{n})U_{n-1} = I - P^{\prime}(x_{n-1})U_{n-1}  - \left(P^{\prime}(x_{n}) -  P^{\prime}(x_{n-1}) \right) U_{n-1}$$
, we obtain:
$$\norm{I - P^{\prime}(x_{n})U_{n-1}} \le \left(A_{n-1}^{2} + \alpha_{n}c_{n-1}\right)
h^{2^{n-1}} = A_{n}h^{2^{n-1}}\, .$$
, from which it emanates that:
$$\norm{U_{n}} \le \norm{U_{n-1}} + \norm{U_{n-1}} \cdot \norm{I - P^{\prime}(x_{n})U_{n-1}} \le c_{n-1}B \left(1 +  A_{n}h^{2^{n-1}}\right) = c_{n}B \, .$$

Now, let us prove that the sequence $\{x_{n}\}$ is fundamental.
\begin{eqnarray}~\label{fundamental}
\begin{tabular}{rcl}
 $\norm{x_{n+p} - x_{n}}$ & $\le$ & $\norm{x_{n+1} - x_{n}} + \ldots + \norm{x_{n+p} - x_{n+p-1}}$ \\
 $ $                                  & $\le$ & $\sum\limits_{k=n+1}^{n+p} \varepsilon_{k}B\eta
                                               \le \varepsilon_{n+1} \left(1+\frac{\varepsilon_{n+2}}{\varepsilon_{n+1}} + \ldots \right) B\eta \le \frac{\varepsilon_{n+1}}{1 - S}\cdot B\eta$ \\
 $ $ & $\le$ & $ \frac{h(1+h)}{2}\cdot \frac{S^{n-1}}{1 - S} \left(N_{1}h\right)^{2^{n}-2}B\eta $
\end{tabular}
\end{eqnarray}
, and as $X$ is a space of type $B$, we know that there exists
$x^{*} = \lim\limits_{n \rightarrow \infty} x_{n}$

From \eqnref{norm x_k to x_0} and the definition of $\gamma$ we obtain, converting to limit, that $x_{*}$ is in the ball $G$.

Note that $P(x)$ is a continuous operator. Thus, by \eqnref{alpha_k}, \eqnref{beta_k}, \eqnref{epsilon_k} and  \eqnref{norm x_k} we obtain:
\begin{eqnarray*}~\label{P_of_x}
\begin{tabular}{rcl}
 $\norm{P(x^{*})}$ & $=$ & $\lim\limits_{n \rightarrow \infty} \norm{P(x_{n})} \le \lim\limits_{n \rightarrow \infty} \beta_{n}h^{2^{n}-1}\eta$ \\
  $ $ & $\le$ & $\lim\limits_{n \rightarrow \infty} \left(\frac{\varepsilon_{n}q_{n-2}^{2}}{c_{1}} + \frac{1}{2}\varepsilon_{n}^{2}h \right)\eta = 0 \, , $
\end{tabular}
\end{eqnarray*}
, from which we can deduce that $P(x^{*}) = 0$. That is, $x^{*}$ is a solution to this equation.

Converting to limit in inequality \eqnref{fundamental}, we obtain the estimate of \eqnref{inTheorem}. Thus, the theorem is proved.\\
\qed
\end{proofi}

\end{document}